\documentclass[a4paper,11pt]{article}
\usepackage{amsfonts}
\usepackage[numbers,longnamesfirst]{natbib}

\newcommand{\E}{\mathbb{E}}
\newcommand{\erre}{\mathbb{R}}

\newcommand{\cp}[2]{\langle#1,#2\rangle}
\newcommand{\ds}{\displaystyle}
\newcommand{\tr}{\mathop{\mathrm{Tr}}\nolimits}

\newtheorem{prop}{Proposition}
\newtheorem{thm}[prop]{Theorem}
\newtheorem{defi}[prop]{Definition}
\newenvironment{demo}{\par\noindent {\em Proof}.\/}{
\hfill\rule{1.5ex}{1.5ex}\medskip\par\noindent}
\newenvironment{rmk}{\addtocounter{prop}{1}%
\medskip\par\noindent\textbf{Remark \theprop\ }}{\medskip\par}

\renewcommand{\cite}[1]{\citet{#1}}

\title{Stochastic optimal control of delay equations\\
arising in advertising models}

\author{Fausto Gozzi\thanks{Dipartimento di Scienze Economiche e
    Aziendali, Libera Universit\`a Internazionale degli Studi Sociali,
    Viale Pola 12, 00198 Roma (Italy).} \ and
    Carlo Marinelli\thanks{Institut f\"ur Angewandte Mathematik, Universit\"at
    Bonn, Wegelerstr. 6, D-53115 Bonn (Germany).}}

\date{December 6, 2004}

\begin{document}
\maketitle

\begin{abstract}
  We consider a class of optimal control problems of stochastic delay
  differential equations (SDDE) that arise in connection with optimal
  advertising under uncertainty for the introduction of a new product to the
  market, generalizing classical work of \cite{NA}. In particular, we deal
  with controlled SDDE where the delay enters both the state and the
  control. Following ideas of \cite{VK} (which however hold only in the
  deterministic case), we reformulate the problem as an infinite dimensional
  stochastic control problem to which we associate, through the dynamic
  programming principle, a second order Hamilton-Jacobi-Bellman equation. We
  show a verification theorem and we exhibit some simple cases where such
  equation admits an explicit smooth solution, allowing us to construct
  optimal feedback controls.
\end{abstract}

\section{Introduction}
In this paper we consider a class of stochastic optimal control
problems where the state equation is a stochastic delay differential
equations (SDDE). Such problems arise for instance in modeling optimal
advertising under uncertainty for the introduction of a new product to
the market, generalizing classical work of \cite{NA}.
The main novelty in our model is that we deal also with delays in
the control: this is interesting from the applied point of view
and introduces new mathematical difficulties in the problem.

Dynamic models in marketing have a long history, that begins at least
with the seminal papers of \cite{VW} and \cite{NA}. Since then a
considerable amount of work has been devoted to the extension of these
models and to their application to problems of optimal advertising, both
in the monopolistic and the competitive settings, mainly under
deterministic assumptions.
Models with uncertainty have received instead relatively less
attention (see \cite{sethi} for a review of the existing work until
1994, \cite{sethi-VWcomp} for a Vidale and Wolf-like model in the
competitive setting, and e.g. \cite{grosset}, \cite{gw} for recent
work on the case of a monopolistic firm).
Moreover, it has been advocated in the literature (see e.g.
\cite{hartl}, \cite{sethi} and references therein), as it reasonable
to assume, that more realistic dynamic models of goodwill should allow
for lags in the effect of advertisement. Namely, it is natural to
assume that there will be a time lag between advertising expenditure
and the corresponding effect on the goodwill level. More generally, a
further lag structure has been considered, allowing a distribution of
the forgetting time.

In this work we incorporate both lag structures mentioned above.
We formulate a stochastic optimal control problem aimed at
maximizing the goodwill level at a given time horizon $T>0$, while
minimizing the cumulative cost of advertising expenditure until
$T$.

This optimization problem is studied using techniques of
stochastic optimal control in infinite dimension. In particular,
we extend to the stochastic case a representation result, proved
by \cite{VK} in the deterministic setting, that allows to
associate to a controlled SDDE with delays both in the state and in
the control a stochastic differential equation (without delays) in
a suitable Hilbert space. This in turn allows us to associate to
the original control problem for the SDDE an equivalent (infinite
dimensional) control problem for the ``lifted'' stochastic
equation.

We deal with the resulting infinite dimensional optimal control
problem through the dynamic programming approach, i.e. through the
study of the associated Hamilton-Jacobi-Bellman (HJB) equation. The
HJB equation that arise in this case is an infinite dimensional second
order semilinear PDE that does not seem to fall into the ones treated
in the existing literature (see Section \ref{sec:delaycontrol} for
details).
Here we give some preliminary results for this equation. First of all
we consider the particular case (but still interesting from the
applied point of view) when the delay does not enter the control term.
In this case the $L^2$ approach of \cite{gozzi-L2} and the
forward-backward SDE approach of \cite{FT1}, \citep{FT2}, \citep{FT}
apply. We show how to apply the former in Section \ref{sec:ex}).

Moreover we consider the general case of delays in the state and in
the control: since we do not know whether a regular solution exist,
the natural approach would be the one of viscosity solutions. We leave
the treatment of the viscosity solution approach to a subsequent work.
Here we concentrate on the special case where regular solutions exist.
In this case we prove a verification theorem and the existence of
optimal feedbacks. Finally, we show through a simple example (but
nevertheless still interesting in applications) that, possibly in
special cases only, smooth solutions may exist, allowing to prove
verification theorems and the existence of optimal feedback controls.

Some further steps are still needed to build a satisfactory theory:
concerning the viscosity solutions theory, it would be important to
get an existence and uniqueness result and a verification theorem;
concerning regular solutions, other cases where further regularity may
arise should be studied. These issues are part of work in progress.

As references on viscosity solutions of HJB equations in infinite
dimensions, their connection with stochastic control and
applications, we refer to e.g. \cite{lions-infdim1},
\citep{lions-infdim2}, \citep{lions-infdim3}, \cite{swiech-unb},
\cite{GRS}, \cite{GS-zakai}, \cite{GSSstoc}.

Other approaches to optimal control problems for systems described by
SDDEs without infinite dimensional reformulation have been proposed in
the literature: for instance, see \cite{elsanosi} and \cite{larssen}
for a more direct application of the dynamic programming principle
without appealing to infinite dimensional analysis, and \cite{KSlibro}
for the linear-quadratic case. See also \cite{EOS} for some solvable
control problems of optimal harvesting, and \cite{EL} for an
application in financial mathematics.

The paper is organized as follows: in section \ref{sec:formul} we
formulate the optimal advertising problem as an optimal control
problem for a SDDE with delays both in the state and the control.
In section \ref{sec:reformul} we prove a representation result
allowing us to ``lift'' this non-Markovian optimization problem to
an infinite dimensional Markovian control problem. In section
\ref{sec:ex} we study the simpler case of a controlled SDDE with
delays in the state only, for which known results apply. Section
\ref{sec:delaycontrol} deals with the general case of delays in
the state and in the control, for which we only give the
verification result. Finally, in section \ref{sec:example} we
construct a simple example of a controlled SDDE with delay in the
state and in the control, whose corresponding HJB equation admits
a smooth solution, hence there exists an optimal control in feedback
form for the control problem.

\section{The advertising model}\label{sec:formul}
Our general reference model for the dynamics of the stock of
advertising goodwill $y(s)$, $0\leq s\leq T$, of the product to be
launched is given by the following controlled stochastic delay
differential equation (SDDE), where $z$ models the intensity of
advertising spending:
\begin{equation}
\label{eq:SDDE} \left\{\begin{array}{l} dy(s) = \ds \left[a_0 y(s)
+ \int_{-r}^0 a_1(\xi)y(s+\xi)\,d\xi
        + b_0 z(s) +
            \int_{-r}^0b_1(\xi)z(s+\xi)\,d\xi\right]ds \\[10pt]
\ds \qquad\qquad  + \sigma\, dW_0(s), \quad 0\leq s\leq T \\[10pt]
y(0)=\eta^0; \quad y(\xi)=\eta(\xi),\;
z(\xi)=\delta(\xi)\;\;\forall\xi\in[-r,0],
\end{array}\right.
\end{equation}
where the Brownian motion $W_0$ is defined on a filtered
probability space
$(\Omega,\mathcal{F},\mathbb{F}=(\mathcal{F}_s)_{s\geq 0},%
\mathbb{P})$, with $\mathbb{F}$ being the completion of the filtration
generated by $W_0$, and $z$ belongs to
$\mathcal{U}:=L^2_\mathbb{F}([0,T];\erre^+)$, the space of square
integrable non-negative processes adapted to $\mathbb{F}$.
Moreover, we shall assume that the following conditions are verified:
\begin{enumerate}
\item $a_0 \leq 0$;
\item $a_1(\cdot) \in L^2([-r,0];\erre)$;
\item $b_0 \geq 0$;
\item $b_1(\cdot) \in L^2([-r,0];\erre_+)$;
\item $\eta^0 \geq 0$;
\item $\eta(\cdot) \geq 0$, with $\eta(0)=\eta^0$;
\item $\delta(\cdot) \geq 0$.
\end{enumerate}
Here $a_0$ is a constant factor of image deterioration in absence
of advertising, $b_0$ is a constant advertising effectiveness
factor, $a_1(\cdot)$ is the distribution of the forgetting time,
and $b_1(\cdot)$ is the density function of the time lag between
the advertising expenditure $z$ and the corresponding effect on
the goodwill level. Moreover, $\eta^0$ is the level of goodwill at
the beginning of the advertising campaign, $\eta(\cdot)$ and
$\delta(\cdot)$ are the histories of the goodwill level and of the
advertising expenditure, respectively, before time zero (one can
assume $\eta(\cdot)=\delta(\cdot)=0$, for instance).
Note that when $a_1(\cdot)$, $b_1(\cdot)$, and $\sigma$ are
identically zero, equation (\ref{eq:SDDE}) reduces to the classical
model of \cite{NA}.

Finally, we define the objective functional
\begin{equation}
\label{eq:obj-orig} J(z(\cdot)) = \E\left[\varphi_0(y(T))-\int_0^T
h_0(z(s))\,ds \right],
\end{equation}
where $\varphi_0$ is a concave utility function with polynomial
growth at infinity, $h_0$ is a convex cost function which is
superlinear at infinity i.e.$$\lim_{z \to +\infty }\frac{h(z)}{z}
=+\infty,$$ and the dynamics of $y$ is determined by
(\ref{eq:SDDE}).
The problem we will deal with is the maximization of the objective
functional $J$ over all admissible controls in $\mathcal{U}$.
%
\begin{rmk}
  Note that in the general framework of delay equations the functions
  $a_1$ and $b_1$ are measures. Here we do not consider such general
  framework for two reasons: first taking $a_1$ and $b_1$ to be $L^2$
  already captures the applied idea of the problem; second, taking
  $a_1$ and $b_1$ to be measures would introduce some technical
  difficulties that we do not want to treat here. More precisely this
  would create some problems in the infinite dimensional reformulation
  bringing unbounded terms into the state equation. Indeed, if
  $b_1\equiv 0$, the case where $a_1$ is a measure can be easily
  treated by a different standard reformulation. This fact allows us
  to treat the case of point delays in the state with no delays in the
  control. This will be the subject of section 4.
\end{rmk}

\section{Reformulation in Hilbert space} \label{sec:reformul}
Throughout the paper, $X$ will be the Hilbert space defined as
$$ X = \erre \times L^2([-r,0];\erre), $$
with inner product
$$ \langle x,y\rangle = x_0y_0 + \int_{-r}^0 x_1(\xi)y_1(\xi)\,d\xi $$
and norm
$$ |x| = \left( |x_0|^2 + \int_{-r}^0 |x_1(\xi)|^2\,d\xi\right)^{1/2}, $$
where $x_0$ and $x_1(\cdot)$ denote the $\erre$-valued and the
$L^2([-r,0];\erre)$-valued components, respectively, of the generic
element $x$ of $X$.

In this section we shall adapt the approach of \cite{VK} to the
stochastic case to recast the SDDE (\ref{eq:SDDE}) as an abstract
stochastic differential equation on the Hilbert space $X$ and so
to reformulate the optimal control problem.

We start by introducing an operator $A:D(A)\subset X\to X$ as
follows:
\begin{eqnarray*}
A: (x_0,x_1(\xi)) &\mapsto& \Big(a_0x_0 + x_1(0),a_1(\xi)x_0
                           -{dx_1(\xi)\over d\xi}\Big) \quad%
\textrm{a.e.}\ \xi\in [-r,0].
\end{eqnarray*}
The domain of $A$ is given by
$$
D(A) = \left\{ (x_0,x_1(\cdot)) \in X:
       x_1(\cdot)\in W^{1,2}([-r,0];\erre),\; x_1(-r)=0\right\}.
$$
The operator $A$ is the adjoint of the operator
$A^*:D(A^*)\subset X \to X$ defined as
\begin{equation}
\label{eq:Astar}
A^*: (x_0,x_1(\cdot)) \mapsto \Big(a_0x_0 + \int_{-r}^0 a_1(\xi)x_1(\xi)\,d\xi,
x'_1(\cdot)\Big)
\end{equation}
on
$$ D(A^*) = \left\{(x_0,x_1(\cdot))\in \erre \times W^{1,2}([-r,0];\erre):
x_0=x_1(0) \right\}. $$ It is well known that $A^*$ is the
infinitesimal generator of a strongly continuous semigroup (see,
e.g., \cite{choj78} or
\cite{DZ96}), therefore the same is true for $A$.\\
We also need to define the bounded linear control operator $B:U \to X$ as
\begin{equation}
B: u \mapsto \Big(b_0u,b_1(\cdot)u\Big),
\end{equation}
where $U:=\erre_+$. Sometimes we shall identify the operator $B$ with
the element $(b_0,b_1(\cdot))\in X$.

We introduce now an abstract stochastic differential equation on the
Hilbert space $X$ that is equivalent, in the sense made precise by
Proposition \ref{prop:equiv}, to the SDDE (\ref{eq:SDDE}):
\begin{equation}
\label{eq:abstract}
\left\{\begin{array}{l}
dY(s) = (AY(s)+Bz(s))\,ds + G\,dW_0(s) \\[8pt]
Y(0) = x \in X,
\end{array}\right.
\end{equation}
where $G:\erre \to X$ is defined by
$$
G: x_0 \to (\sigma x_0, 0),
$$
and $z(\cdot) \in \mathcal{U}$. Using theorems 5.4 and 5.9 in
\cite{DZ92}, we have that equation (\ref{eq:abstract}) admits a unique
weak solution (in the sense of \citep{DZ92}, p.~119) with continuous paths
given by the variation of constants formula:
$$
Y(s) = e^{sA}x + \int_0^s e^{(s-\tau)A}Bz(\tau)\,d\tau + \int_0^s
e^{(s-\tau)A}G\,dW_0(\tau).
$$
In order to state equivalence results between the SDDE (\ref{eq:SDDE})
and the abstract evolution equation (\ref{eq:abstract}), we need to
introduce the operator $M: X \times L^2([-r,0];\erre) \to X$ defined
as follows:
$$ M:(x_0,x_1(\cdot),v(\cdot)) \mapsto (x_0,m(\cdot)), $$
where
$$
m(\xi) := \int_{-r}^\xi a_1(\zeta) x_1(\zeta-\xi)\,d\zeta
+ \int_{-r}^\xi b_1(\zeta) v(\zeta-\xi)\,d\zeta.
$$
The following result is a generalization of Theorems 5.1 and 5.2
of \cite{VK}.
\begin{prop}
\label{prop:equiv} Let $Y(\cdot)$ be the weak solution of the abstract
evolution equation (\ref{eq:abstract}) with arbitrary initial
datum $x\in X$ and control $z(\cdot) \in {\cal U}$. Then, for $t\geq r$,
one has, $\mathbb{P}$-a.s.,
$$ Y(t) = M(Y_0(t),Y_0(t+\cdot),z(t+\cdot)). $$
Moreover, let $\{y(t),\; t\geq -r\}$ be a continuous solution of the
stochastic delay differential equation (\ref{eq:SDDE}), and $Y(\cdot)$
be the weak solution of the abstract evolution equation
(\ref{eq:abstract}) with initial condition
$$ x = M(\eta^0,\eta(\cdot),\delta(\cdot)). $$
Then, for $t\geq 0$, one has, $\mathbb{P}$-a.s.,
$$ 
Y(t) = M(y(t),y(t+\cdot),z(t+\cdot)),
$$ 
hence $y(t)=Y_0(t)$, $\mathbb{P}$-a.s., for all $t\geq 0$.
\end{prop}
\begin{demo}
  Let $x=(x_0,x_1)\in D(A)$ (for general $x$ the same result will
  follow by standard density arguments -- see e.g. \cite{VK}).
Equation (\ref{eq:abstract}) can be written as
\begin{equation}
\left\{
\begin{array}{l}
\ds dY_0(t) = \Big(a_0Y_0(t)+Y_1(t)(0) + b_0z(t)\Big)\,dt
    + \sigma\,dW_0(t)\\[8pt]
\ds dY_1(t)(\xi) = \Big(a_1(\xi)Y_0(t) - {d \over d\xi}Y_1(t)(\xi)
    + b_1(\xi)z(t)\Big)\,dt,\\[8pt]
\ds Y_0(0) = x_0, \quad Y_1(0)(\cdot) = x_1(\cdot),
\end{array}\right.
\end{equation}
therefore, $\mathbb{P}$-a.s.,
\begin{eqnarray}
Y_0(t) &=& e^{ta_0}x_0 + \int_0^t e^{(t-s)a_0}Y_1(s)(0)\,ds +
\int_0^t e^{(t-s)a_0}b_0z(s)\,ds + \int_0^t
e^{(t-s)a_0}\sigma\,dW_0(s)\nonumber \\
\label{eq:Y1} 
Y_1(t)(\xi) &=& [\Phi(t)x_1](\xi) + \int_0^t
[\Phi(t-s)a_1(\cdot)Y_0(s)](\xi)\,ds + \int_0^t
[\Phi(t-s)b_1(\cdot)z(s)](\xi)\,ds,
\end{eqnarray}
where $\Phi(t)$ is the semigroup of truncated right shifts defined as
$$
[\Phi(t)f(\cdot)](\xi) = \left\{
\begin{array}{ll}
f(\xi-t), & -r \leq \xi-t \leq 0,\\
0, & \textrm{otherwise}
\end{array}\right.
$$
for all $f\in L^2$. Then (\ref{eq:Y1}) for $t\geq r$ can be
rewritten, using the definition of $\Phi$ and recalling that
$x_1(-r)=0$, since $x\in D(A)$, as ($\mathbb{P}$-a.s.)
\begin{equation}
\label{eq:seconda}
Y_1(t)(\xi) = \int_{-r}^\xi a_1(\alpha)Y_0(t+\alpha-\xi)\,d\alpha
+ \int_{-r}^\xi b_1(\alpha)z(t+\alpha-\xi)\,d\alpha,
\end{equation}
which is equivalent to
$$
Y(t) = M(Y_0(t),Y_0(t+\cdot),z(t+\cdot)),
$$
as claimed. 

Let us now prove the second claim: the $L^2$-valued component of
the weak solution of the evolution equation (\ref{eq:abstract})
with initial data $Y(0)=x=M(\eta^0,\eta(\cdot),\delta(\cdot)).$
satisfies, for $\xi-t\in[-r,0]$, $\xi\in[-r,0]$, $t \geq 0$,
$$
Y_1(0)(\xi-t) = \int_{-r}^{\xi-t} a_1(\alpha)\eta(t+\alpha-\xi)\,d\alpha
+ \int_{-r}^{\xi-t} b_1(\alpha)\delta(t+\alpha-\xi)\,d\alpha,
$$
$\mathbb{P}$-a.s., as follows from (\ref{eq:Y1}). We assume here
$\eta(0)=\eta^0$, without loss of generality (the general case follows
by density arguments, as in \citep{VK}). Again by (\ref{eq:Y1}) and
some calculations we obtain, $\mathbb{P}$-a.s.,
$$
Y_1(t)(\xi) = \int_{-r}^\xi a_1(\alpha)\tilde{Y}_0(t+\alpha-\xi)\,d\alpha +
\int_{-r}^\xi b_1(\alpha)z(t+\alpha-\xi)\,d\alpha,
$$
where
$$
\tilde{Y}_0(\xi) = \left\{
\begin{array}{ll}
\eta(\xi), & \xi \in [-r,0],\\ 
Y_0(\xi), & \xi \geq 0.
\end{array}\right.
$$
Observe that the definition of $\tilde{Y}$ is well posed because
of the assumption $\eta(0)=\eta^0$, and because $\eta^0=Y_0(0)$ by
the definition of the operator $M$.
In order to finish, we need to prove that
$Y_0(\cdot)$ satisfies the same integral equation (in the mild sense) as the solution $y(\cdot)$ to the SDDE (\ref{eq:SDDE}), i.e. that the following
holds for all $t\geq 0$:
$$
\int_0^t e^{(t-s)a_0}Y_1(s)(0)\,ds =
\int_0^t e^{(t-s)a_0} \left[
\int_{-r}^0 a_1(\xi)Y_0(s+\xi)\,d\xi + \int_{-r}^0 b_1(\xi)z(s+\xi)\,d\xi
\right]\,ds.
$$
But this follows immediately from (\ref{eq:seconda}) with $\xi=0$:
$$
Y_1(t)(0) =
\int_{-r}^0 a_1(\xi)Y_0(t+\xi)\,d\xi + \int_{-r}^0 b_1(\xi)z(t+\xi)\,d\xi,
$$
which proves the claim. The fact that $y(t)=Y_0(t)$,
$\mathbb{P}$-a.s., for all $t\geq 0$, easily follows.
\end{demo}

Using Proposition \ref{prop:equiv}, we can now give a Hilbert space
formulation of our problem. Since we want to use the dynamic
programming approach, from now on we let the initial time vary,
calling it $t$ with $0\le t\le T$.

The state space is $X = \erre \times L^2([-r,0];\erre) $, the
control space is $U=\erre_+$, the control strategy is $z(\cdot)\in
{\cal U}$.
The state equation is (\ref{eq:abstract}) with initial condition
at $t$, i.e.
\begin{equation}
\label{eq:abstractbis}
\left\{\begin{array}{l}
dY(s) = (AY(s)+Bz(s))\,ds + G\,dW_0(s) \\[8pt]
Y(t) = x \in X,
\end{array}\right.
\end{equation}
and its unique weak solution, given the initial data $(t,x)$ and
the control strategy $z(\cdot)$, will be denoted by
$Y(\cdot;t,x,z(\cdot))$.
The objective functional is
\begin{equation}
\label{ex1jbis} J(t,x;z(\cdot)) =
\E^{t,x}\left[\varphi(Y(T,t,x;z(\cdot)))+\int_t^T h(z(s))\,ds
\right].,
\end{equation}
where the functions $h:U\to\erre$ and $\varphi:X\to\erre$ are defined as
\begin{eqnarray*}
h(z) &=& -h_0(z) \\
\varphi(x_0,x_1(\cdot)) &=& \varphi_0(x_0).
\end{eqnarray*}
The problem is to maximize the objective function
$J(t,y;z(\cdot))$ over all $z(\cdot)\in \mathcal{U}$.
We also define the value function $V$ for this problem as
$$
V(t,x) = \inf_{z(\cdot)\in\mathcal{U}} J(t,x;z(\cdot)). 
$$
Moreover, we shall say that $z^*(\cdot) \in \mathcal{U}$ is an optimal
control if it is such that
$$ 
V(t,x) = J(t,x;z^*(\cdot)). 
$$
Following the dynamic programming approach we would like first to
characterize the value function as the unique solution (in a
suitable sense) of the following Hamilton-Jacobi-Bellman equation
\begin{equation}
\label{eq:HJB0} \left\{\begin{array}{l}
\ds v_t + {1\over 2}\tr (Qv_{xx}) + \cp{Ax}{v_x} + H_0(v_x) = 0 \\[8pt]
v(T) = \varphi,
\end{array}\right.
\end{equation}
where $Q=G^*G$ and
$$
H_0(p) = \sup_{z\in U} (\cp{Bz}{p}+h(z)).
$$
Moreover we would like to find a sufficient condition for
optimality given in terms of $V$ (a so-called verification theorem)
and a feedback formula for the optimal control $z^*$.

\section{The case with no delay in the control}
\label{sec:ex}
In a model for the dynamics of goodwill with distributed
forgetting factor, but without lags in the effect of advertising
expenditure, i.e. with $b_1(\cdot)=0$ in (\ref{eq:SDDE}), it is
possible to apply both the approach of Hamilton-Jacobi-Bellman
equations in $L^2$ spaces developed by \cite{gozzi-L2}, and the
backward SDE approach of \cite{FT}. We follow here the first
approach, showing that both the value function and the optimal
advertising policy can be characterized in terms of the solution
of a Hamilton-Jacobi-Bellman equation in infinite dimension. In
fact we treat a somewhat different case assuming that the
distribution of the forgetting factor is concentrated on a point.
This can be treated by a different standard reformulation and with
simpler computations and interpretations.
In particular, we consider the case where the goodwill evolves
according to the following equation:
\begin{equation}
\label{eq:ex1} \left\{\begin{array}{l} dy(s) = \ds [a_0 y(s) +
a_1y(s-r)
        + b_0 z_0(s)]\,ds + \sigma\, dW_0(s), \quad 0\leq s\leq T \\[10pt]
y(0)=\eta^0; \quad y(\xi)=\eta(\xi) \;\; \forall\xi\in[-r,0].
\end{array}\right.
\end{equation}
By the representation theorems for solutions of stochastic delay
equations of \cite{choj78}, one can associate to (\ref{eq:ex1}) an
evolution equation on the Hilbert space
$X$ of the type
\begin{equation}
\label{eq:ex1a}
\left\{\begin{array}{l}
dY(s) = (AY(s) + \sqrt{Q}z(s))\,ds + \sqrt{Q}\,dW(s), \\[10pt]
Y(0)=y,
\end{array}\right.
\end{equation}
where $A:D(A)\subset X \to X$ is defined as
$$ A:(x_0,x_1(\cdot)) \mapsto (a_0x_0+a_1x_1(-r),x_1'(\cdot)) $$
on its domain $$D(A) = \left\{(x_0,x_1(\cdot))\in (\erre \times
W^{1,2}([-r,0];\erre): x_0=x_1(0) \right\};$$ moreover
$z=(\sigma^{-1}b_0z_0,z_1(\cdot)) \in \erre_+\times
L^2([-r,0];\erre)$, with $z_1(\cdot)$ a fictitious control; $Q:X
\to X$ is defined as
$$ Q: (x_0,x_1(\cdot)) \mapsto (\sigma^2x_0,0); $$
$W$ is an $X$-valued cylindrical Wiener process with
$W=(W_0,W_1)$, and $W_1$ is a (fictitious) cylindrical Wiener
process taking values in $L^2([-r,0];\erre)$. Finally,
$y=(\eta^0,\eta(\cdot))$.
\begin{rmk}
  Note that the operator $A$ just introduced does not coincide with
  the $A$ introduced in section \ref{sec:reformul}. In fact, $A$ here
  is exactly the $A^*$ defined there . Similarly, the initial datum of
  this section differs from that of section \ref{sec:reformul}. Note
  also that the reformulation carried out in this section does not
  extend to the more general case of delay in the control, explaining
  why the more elaborate construction of the previous section is
  needed.
  
  We also note that the insertion of the fictitious control $z_1$ is
  not necessary here. We do it to keep the control space $U$ equal to
  the state space $X$ so the formulation falls into the results
  contained in \cite{gozzi-L2}. However, it can be easily proved that
  the results of \citep{gozzi-L2} still hold when the weaker condition
  $B(U)\subset Q^{\frac12}(X)$ is satisfied.
\end{rmk}
The operator $A$ is the infinitesimal generator of a strongly
continuous semigroup $\{S(s), s\ge 0\}$ (see again \cite{choj78}).
More precisely, one has
$$ S(s)(x_0,x_1(\cdot)) = (y(s),y(s+\xi)|_{\xi\in[-r,0]}), $$
where $y(\cdot)$ is the solution of the deterministic delay equation
\begin{equation}
\left\{\begin{array}{l}
\ds {dy(s)\over dt} = a_0 y(s) + a_1y(s-r), \quad 0\leq s\leq T \\[10pt]
y(0)=x_0; \quad y(\xi)=x_1(\xi) \;\; \forall\xi\in[-r,0].
\end{array}\right.
\end{equation}
Moreover, the $X$-valued mild solution
$Y(\cdot)=(Y_0(\cdot),Y_1(\cdot))$ of (\ref{eq:ex1a}) is such that
$Y_0(\cdot)$ solves the original stochastic delay equation
(\ref{eq:ex1}).

As in section \ref{sec:reformul}, we now consider an associated
stochastic control problem letting the initial time $t$ vary in
$[0,T]$. The state equation is (\ref{eq:ex1a}) with initial
condition at $t$, i.e.
\begin{equation}
\label{eq:ex1abis} \left\{\begin{array}{l}
dY(s) = (AY(s) + \sqrt{Q}z(s))\,ds + \sqrt{Q}\,dW(s), \\[10pt]
Y(t)=x,
\end{array}\right.
\end{equation}
and the objective function is
$$
J(t,x;z_0(\cdot)) = \E^{t,x}\left[
\varphi_0(y(T)) - \int_t^T h_0(z_0(s))\,ds\right],
$$
with $y(\cdot)$ obeying the SDDE (\ref{eq:ex1}). Defining
\begin{eqnarray*}
h(z_0,z_1(\cdot)) &=& -h_0(z_0) \\
\varphi(x_0,x_1(\cdot)) &=& \varphi_0(x_0),
\end{eqnarray*}
thanks to the above mentioned equivalence between the SDDE
(\ref{eq:ex1}) and the abstract SDE (\ref{eq:ex1a}), we are led to the
equivalent problem of maximizing the objective function
\begin{equation}
\label{ex1j} J(t,x;z(\cdot)) =
\E^{t,x}\left[\varphi(Y(T))+\int_t^T h(z(s))\,ds \right],
\end{equation}
with $Y$ subject to (\ref{eq:ex1a}).\\
Before doing so, however, following \citep{gozzi-L2}, we need to
assume conditions on the coefficients of (\ref{eq:ex1}) such that
the uncontrolled version of (\ref{eq:ex1a}), i.e.
\begin{equation}
\label{eq:ex1auc}
\left\{\begin{array}{l}
dZ(t) = AZ(t)\,dt + \sqrt{Q}\,dW(t), \\[10pt]
Z(0)=x,
\end{array}\right.
\end{equation}
admits an invariant measure. It is known (see \cite{DZ96}) that
\begin{equation}
\label{ex1-cond}
a_0 < 1, \quad a_0 < -a_1 < \sqrt{\gamma^2+a_0^2},
\end{equation}
where $\gamma\,\mathrm{coth}\,\gamma = a_0$, $\gamma\in]0,\pi[$, ensures that
(\ref{eq:ex1auc}) admits a unique non-degenerate invariant measure
$\mu$.
\begin{rmk}
  The deterioration factor $a_0$ is always assumed to be negative,
  hence the first condition in (\ref{ex1-cond}) is not a real
  constraint. In general, however, it is not clear what sign $a_1$
  should have. If $a_1$ is also negative, i.e. it can again be
  interpreted as a deterioration factor, condition (\ref{ex1-cond})
  says that $a_1$ cannot be ``much more negative'' than $a_0$. On the
  other hand, if $a_1$ is positive, then the second condition in
  (\ref{ex1-cond}) implies that the improvement effect as measured by
  $a_1$ cannot exceed the deterioration effect as measured by $|a_0|$.
\end{rmk}
Let us now define the Hamiltonian $H_0:X \to \erre$ as
$$ H_0(p) = \sup_{z\in U} \Big(\cp{\sqrt{Q}z}{p}_X + h(z)\Big) $$
and write the Hamilton-Jacobi-Bellman equation associated to the control
problem (\ref{ex1j}):
\begin{equation}
\label{ex1hjb}
\left\{\begin{array}{ll}
\ds v_t + {1\over 2}\tr(Qv_{xx}) + \cp{Ax}{v_x} + H_0(v_x) = 0 \\[10pt]
v(T,x) = \varphi(x).
\end{array}\right.
\end{equation}
If the Hamiltonian $H_0$ is Lipschitz (which follows from the
hypothesis on $h_0$), $\varphi\in L^2(X,\mu)$ (which follows from
the hypothesis on $\varphi$), and the operator $A$ satisfies
(\ref{ex1-cond}), then (\ref{ex1hjb}) admits a unique mild
solution $v$ in the space $L^2(0,T; W^{1,2}_Q(X,\mu))$, as follows
from Theorem 3.7 of \citep{gozzi-L2}. Moreover, Theorem 5.7 of
ibid. guarantees that $v$ coincides ($\mu$-a.e.) with the value
function
$$
V(t,x) = \inf_{z\in\mathcal{U}}
\E^{t,x}\left[\varphi(Y(T))+\int_t^T h(z(s))\,ds \right]
$$
(by $z\in\mathcal{U}$ we mean, with a slight abuse of notation,
$z_0\in\mathcal{U}$), and that there exists a unique optimal control
$z^*$, i.e.
$$ V(t,x) = J(t,x;z^*(\cdot)) $$
with
\begin{equation}
\label{ex1-zopt}
 z^*(t) = DH(\widetilde{D}_Q v(t,Y^*(t))),
\end{equation}
and $Y^*$ is the solution of the closed-loop equation
\begin{equation}
\label{ex1-cl} \left\{\begin{array}{ll} \ds dY(s) = [AY(s) +
\sqrt{Q}DH(\widetilde{D}_Qv(s,Y(s))]\,ds
            + \sqrt{Q}\,dW(t) \\[10pt]
Y(t)=x.
\end{array}\right.
\end{equation}
The gradient operator $\widetilde{D}_Q$ is, roughly speaking, a
``weakly closable'' extension of the Malliavin derivative $Q^{1/2}D$
acting on $W_Q^{1,2}(X,\mu)$. For the exact definition and
construction of $\widetilde{D}_Q$ we refer the reader to
\citep{gozzi-L2}.
\begin{rmk}
  The HJB equation (\ref{ex1hjb}) is ``genuinely'' infinite
  dimensional, i.e. it reduces to a finite dimensional one only in
  very special cases. For example, by the results of \cite{lari},
  (\ref{ex1hjb}) reduces to a finite dimensional PDE if and only if
  $a_0=-a_1$. However, under this assumption, we cannot guarantee the
  existence of a non-degenerate invariant measure for the
  Ornstein-Uhlenbeck semigroup associated to (\ref{eq:ex1auc}).  Even
  more extreme would be the situation of distributed forgetting time:
  in this case the HJB is finite dimensional only if the term
  accounting for distributed forgetting vanishes altogether!
\end{rmk}

\section{Delays in the state and in the control}
\label{sec:delaycontrol}
We now consider the case when also delays in the control
are present. The optimal control problem is the one described in
section 2 with $a_1(\cdot)\neq 0$, $b_1(\cdot) \neq 0$.

The Hamilton-Jacobi-Bellman equation associated to the problem is
\begin{equation}
\label{eq:HJB} \left\{\begin{array}{l}
\ds v_t + {1\over 2}\tr (Qv_{xx}) + \cp{Ax}{v_x} + H_0(v_x) = 0 \\[8pt]
v(T) = \varphi,
\end{array}\right.
\end{equation}
where $Q=G^*G$ and
$$
H_0(p) = \sup_{z\in U} (\cp{Bz}{p}+h(z)).
$$

Unfortunately it is not possible, in general, to obtain regular
solutions of the HJB equation (\ref{eq:HJB}) using the existing
theory based on perturbations of Ornstein-Uhlenbeck semigroups
(see e.g. \cite{BDP}, \cite{DZ02} and \cite{gozzi95},
\citep{gozzi96}).  The main problem is the lack of regularity
properties of a suitable Ornstein-Uhlenbeck semigroup associated
to the problem: in particular, the associated gradient operator is
not closable and the semigroup is not strongly Feller (see
\cite{gozzi-L2} and \cite{gozzi-closab}).

As mentioned in section \ref{sec:ex}, if there is no lag in the effect of
advertisement spending on the goodwill, i.e. if $b_1(\cdot)=0$,
then both the $L^2$ approach of \cite{gozzi-L2} and the backward
SDE approach of \cite{FT} can be applied, obtaining a
characterization of the value function and of the optimal
advertising policy in terms of the (unique) solution to
(\ref{eq:HJB}). Both approaches, however, fail in this more
general case, since they require, roughly speaking, that the image
of the control is included in the image of the noise, i.e. that
$B(U) \subseteq G(X)$, which is clearly not true.

The only approach that seems to work in the general case of delays
in the state and in the control is, to the best of our knowledge,
the framework of viscosity solutions (see \cite{lions-infdim1},
\citep{lions-infdim2}, \citep{lions-infdim3}).  However, while
this approach gives a characterization of the value function in
terms of the unique (viscosity) solution of the
Hamilton-Jacobi-Bellman equation (\ref{eq:HJB}), this solution is
only guaranteed to be continuous, hence one can construct from it
an optimal control only in a rather weak form, through the so
called viscosity verification theorems (see \cite{GSZ}). The study
of this problem in the framework of viscosity solutions is the
subject of a forthcoming publication.

Here we want to prove a verification theorem in the case when
regular solutions of the HJB equation are available.

\begin{defi}
\label{solclass}
A function $v$ is said to be 
\begin{itemize}
\item a \emph{classical solution} of the HJB equation (\ref{eq:HJB}) if $v\in
C^{1,2}([0,T]\times X)$ and $v$ satisfies (\ref{eq:HJB})
pointwise;
\item an \emph{integral solution} if $v\in C^{0,2}([0,T]\times X)$,
and moreover for $t \in [0,T]$ and $x\in D(A)$, one has
\begin{equation}
\label{eq:HJBint}
\varphi (x) -v(t,x) + \int_t^T \left[{1\over 2}\tr (Qv_{xx}(s,x))
+ \cp{Ax}{v_x(s,x)} + H_0(v_x(s,x))\right]ds = 0
\end{equation}
\end{itemize}
\end{defi}

\begin{thm}[Verification Theorem]
\label{th:VT}
Let $v$ be an integral solution of the HJB (\ref{eq:HJB}) and let $V$
be the value function of the optimal control problem. Then
\begin{itemize}
\item[(i)]  $v\geq V$ on $ \left[ 0,T\right]\times X$.
\item[(ii)]  If a control $z\in \mathcal{U} $ is such that, at
starting point $(t,x)$,
$$
H_0(v_x(s,Y(s))) = \sup_{z\in U} \{\cp{Bz}{v_x(s,Y(s))}+h(z)\}
=\cp{Bz(s)}{v_x(s,Y(s))}+h(z(s))
$$
for almost every $s\in \left[ t,T\right]$, ${\mathbb P}$-a.e.,
then this control is optimal and $v(t,x) =V(t,x)$.
\item[(iii)]  If we know a priori that $V=v$, then (ii) is a necessary (and
sufficient) condition of optimality.
\end{itemize}
\end{thm}

Although there is a standard way to prove such results, this
version of the verification theorem is not contained in the existing
literature.

We give an idea of the method by sketching the proof in the case
of bounded $A$. The case of unbounded $A$ can be treated by
approximating $A$ with its Yosida approximations $A_{n}$,  and
then passing to the limit as $n\rightarrow +\infty $ (see e.g. \cite{BDP}).

\medskip

\begin{demo}
Let $A$ and $B$ be bounded operators. The core of
the job is to prove that, for every $\left( t,x\right) \in \left[
0,T\right] \times X$ and any $z\in \mathcal{U}$ the following
fundamental identity holds
\begin{equation}\label{eq:fundid}
v\left( t,x\right)=J\left( t,x;z\left( \cdot \right) \right)
 +\int_{t}^{T}
\left[H_0(v_x(s,Y(s))) -\cp{Bz(s)}{v_x(s,Y(s))}-h(z(s))\right]ds,
\end{equation}
where $Y(s):=Y(s;t,x;z(\cdot))$. Once this is proved, the three
claims of the theorem follow as straightforward consequences of the
definitions of the Hamiltonian $H_0$, of value function and of optimal
strategy.

Let us then prove (\ref{eq:fundid}). We first approximate $v$ by a
sequence of smooth function $v_n$ that solve suitable
approximating equations and that are such that
$$
v_n \longrightarrow v, \qquad v_{nx} \longrightarrow v_x.
$$
This is possible e.g. using the same ideas of
\cite{gozzi-L2}, Section 4. Then for a.e. $s\in \left( \left[
    t,T\right] \cap \mathbb{R}\right)$, one may show that
\begin{eqnarray}
\frac{d}{ds} v_n(s,Y(s)) &=& v_{nt}(s,Y(s)) + \cp{y'(s)}{v_{nx}(s,Y(s))} \\
&=& v_{nt}(s,Y(s)) + \cp{AY(s)+Bz(s)}{v_{nx}(s,Y(s))}.
\label{eq:dertot}
\end{eqnarray}
Since $v_n$ is a classical solution of a suitable approximating HJB
equation (see again \citep{gozzi-L2}), we have
$$
v_{nt}(s,Y(s)) + \cp{AY(s)}{v_{ns}(s,Y(s))} =
-H_0(B^\ast v_{nx}(s,Y(s))) - g_n(s,Y(s)),
$$
where $g_n$ is a term appearing in the approximating HJB such that
$g_n \to 0$ as $n \to \infty$ (see again \citep{gozzi-L2}).
Substituting in (\ref{eq:dertot}) and then adding and subtracting
$h\left(s,z\left( s\right) \right)$, we obtain
\begin{eqnarray*}
\frac{d}{ds}v_n(s,Y(s)) &=& \cp{z(s)}{B^\ast v_{nx}(s,Y(s))}
-H_0(B^{\ast }v_{nx}(s,Y(s))) - g_n(s,Y(s))\\
&=& \cp{z(s)}{B^{\ast}v_{nx}(s,Y(s))} + h(z(s)) \\
&& -H_0(B^{\ast }v_{nx}(s,Y(s))) - g_n(s,Y(s)) - h(z(s)).
\end{eqnarray*}
Integrating between $t$ and $T$ we get
\begin{eqnarray*}
&&v_n\left( T,Y\left( T\right) \right) -v_n\left( t,x\right)
+\int_{t}^{T}\left[
g_n\left( s,Y\left( s\right) \right) +h(z(s)) %
\right] ds \\
&=& \int_t^T \left[ \cp{z(s)}{B^*(v_{nx}(s,Y(s)))} + h(z(s))
- H_0(B^*v_{nx}(s,Y(s))) \right]\,ds,
\end{eqnarray*}
which gives the desired equality (\ref{eq:fundid}) after passing to
the limit as $n \to \infty$.
\end{demo}
About feedback maps, the following theorem holds true.
\begin{thm}
\label{th:feedback}
Under the same hypotheses of theorem \ref{th:VT}, assume that there
exists a measurable map $\Gamma:\left(\left[ 0,T\right]\cap
  \mathbb{R}\right) \times X\rightarrow U$, such that $\Gamma\left(
  t,p\right) $ is a maximum point of the map $z\mapsto\langle
z,B^*p\rangle+h(z)$ for given $\left( t,p\right) \in \left(\left[
    0,T\right]\cap \mathbb{R}\right) \times X$. Assume also that for
every $\left( t,x\right) $ there exists a (mild) solution $Y^{\ast
}\left( \cdot \right) $ of the closed loop equation
$$
\left\{
\begin{array}{ll}
\ds dY^{\ast}(\tau)=AY^*(\tau ) + B\Gamma(\tau,v_{x}(\tau,Y^{\ast }(\tau)))
+ G\,dW(t),
& \tau \in \left[ t,T\right]\cap \mathbb{R}; \\[8pt]
Y^*(t)=x, & x\in X,
\end{array}
\right.
$$
such that $\tau \mapsto \Gamma(\tau,v_{x}(\tau,Y^*(\tau)))$ belongs
to $\mathcal{U}$. Then the couple $(Y^*(\cdot),z^*(\cdot))$ is
optimal whenever the control strategy satisfies the feedback relation
$$
z^{\ast }\left( \tau \right) =\Gamma\left( \tau, v_{x}\left( \tau
,Y^{\ast }\left( \tau \right) \right) \right).
$$
\end{thm}

\begin{demo} This is a straightforward consequence of (ii) in
the Verification Theorem, as the control generated by the feedback
relation satisfies (ii).
\end{demo}


\section{An example with explicit solution}
\label{sec:example}
Here we restrict ourselves to some less general specification of
the objective function, but still meaningful, for which the HJB
equation admits a smooth solution, and therefore the value
function and the optimal control can be completely characterized.

Let us assume
$$ h(z)=-\beta z_0^2 $$
and
$$ \varphi(x) = \gamma x_0, $$
with $\beta$, $\gamma>0$.
Then we have
$$ H_{CV}(p,z) = \cp{Bz}{p}+h(z) = \cp{B}{p}z - \beta z^2, $$
and
$$
H_0(p) = \sup_{z\in U} H_{CV}(z,p) =%
\left\{\begin{array}{ll}
\ds {\cp{B}{p}^2 \over 4\beta}, & \cp{B}{p} \geq 0 \\[10pt]
\ds 0, & \cp{B}{p} < 0,
\end{array}\right.
$$
or equivalently, in more compact notation,
$$ H_0(p) = {(\cp{B}{p}^+)^2 \over 4\beta}. $$
We conjecture a solution of the HJB equation (\ref{eq:HJB}) of the form
$$ 
v(t,x) = \cp{w(t)}{x} + c(t), \qquad t \in [0,T], \; x \in X,
$$
where $w(\cdot):[0,T]\to X$ and $c(\cdot):[0,T]\to\erre$
are given functions whose properties we will study below. Hence
for $t \in [0,T] $ and $x \in X$ we have, assuming that all
objects are well defined,
\begin{eqnarray*}
v_t(t,x) &=& \cp{w'(t)}{x} + c'(t) \\
v_x(t,x) &=& w(t) \\
v_{xx} &=& 0,
\end{eqnarray*}
and
\begin{equation}\label{eq:HJBf}
\left\{\begin{array}{ll}
\ds \cp{w'(t)}{x} + c'(t) + \cp{Ax}{w} +
{(\cp{B}{w(t)}^+)^2 \over 4\beta} = 0, & t \in [0,T), \; x \in X,\\[10pt]
\ds \cp{w(T)}{x} + c(T) = \gamma x_0, & x \in X.
\end{array}\right.
\end{equation}
It is clear that, if $A^*w(t)$ is well defined for all $t\in [0,T]$,
(\ref{eq:HJBf}) is equivalent to
\begin{equation}
\label{eq:HJBf1}
\left\{\begin{array}{ll}
\cp{w'(t)}{x} + \cp{A^*w(t)}{x} = 0, & t \in [0,T[ \\[8pt]
w(T)=(\gamma, 0)
\end{array}\right.
\end{equation}
and
\begin{equation}
\label{eq:HJBf2}
\left\{\begin{array}{ll}
\ds c'(t) + {(\cp{B}{w(t)}^+)^2 \over 4\beta} = 0, & t \in [0,T[ \\[8pt]
c(T)=0.
\end{array}\right.
\end{equation}
Since (\ref{eq:HJBf1}) must hold for all $x$, then it implies
$$
\left\{\begin{array}{ll}
w'(t) + A^*w(t) = 0, & t \in [0,T[ \\[8pt]
w(T)=(\gamma, 0). 
\end{array}\right.
$$
Recalling (\ref{eq:Astar}), we obtain that (\ref{eq:HJBf1}) is equivalent to
\begin{equation}\label{eq:HJBf3}
\left\{
\begin{array}{ll}
\ds w_0'(t) + a_0w_0(t) + \int_{-r}^0 a_1(\xi)w_1(t,\xi)\,d\xi = 0, &
t \in [0,T[ \\[8pt]
w_0(T)=\gamma
\end{array}\right.
\end{equation}
and
\begin{equation}\label{eq:HJBf4}
\left\{
\begin{array}{ll}
\ds {\partial w_1 \over \partial t}(t,\xi) +%
{\partial w_1 \over \partial \xi}(t,\xi) = 0, &
t \in [0,T[, \; \xi \in [-r,0[ \\[8pt]
w_1(T,\xi) = 0, & \xi \in [-r,0[ \\[8pt]
w_1(t,0) = w_0(t), & t \in [0,T].
\end{array}\right.
\end{equation}
The solution of (\ref{eq:HJBf4}) is given by
\begin{equation}
w_1(t,\xi) = w_0(t-\xi)\mathbb{I}_{\{t-\xi\in[0,T]\}},
\end{equation}
from which one can solve equation (\ref{eq:HJBf3}), obtaining
$w_0(\cdot)$. Note that, unfortunately, the function $w$ is never in
$D(A^*)$, except for the case when it is $0$ everywhere. However this
does not exclude that the candidate $v(t,x)=\cp{w(t)}{x}+c(t)$ solves
the HJB equation (\ref{eq:HJB}) in some sense. Indeed it is an
integral solution of (\ref{eq:HJB}) in the sense of definition
(\ref{solclass}), as it follows from the above calculations. Note that
$v\in C([0,T]\times X)$ and that it is twice differentiable in the $x$
variable, i.e. it satisfies the hypotheses of the verification theorem
\ref{th:VT}. Since a maximizer of the current-value Hamiltonian is
given by $z^*=\cp{B}{p}^+/(2\beta)$, then it is immediately seen that
the control
$$
z^*(t) = \frac{\cp{B}{v_x(t)}^+}{2\beta} =
\frac{\cp{B}{w(t)}^+}{2\beta}, \qquad t \in [0,T],
$$
which does not depend on $Y^*(t)$, is admissible, hence (ii) of the
verification theorem \ref{th:VT} holds, and $z^*(\cdot)$ is optimal.


\def\cprime{$'$}

\end{document}